\documentclass[10pt,a4paper]{article}
\usepackage{amsmath}
\usepackage{latexsym}
\usepackage{amssymb}
\usepackage{amstext}
\usepackage{amsfonts}
\usepackage{amsthm}
\usepackage{fancyhdr}
\usepackage{indentfirst}
\usepackage{hyperref}
\usepackage{graphicx}
\usepackage{epsfig}
\usepackage{pgf}
\usepackage{mathrsfs}
\usepackage{listings}
\usepackage{cases}
\usepackage{epstopdf}
\usepackage{graphicx}
\usepackage{enumerate}
\usepackage{abstract}
\usepackage{caption}
\date{}

\begin{document}
	\title{Construction of graphs being determined by their generalized $Q$-spectra}
	\author{Gui-Xian Tian$^a$\footnote{Corresponding author. E-mail address: gxtian@zjnu.cn or guixiantian@163.com.}, Jun-Xing Wu$^a$, Shu-Yu Cui$^{b}$, Hui-Lu Sun$^a$ \\%EndAName
		{\small{\it $^a$Department of Mathematics,}}
		{\small{\it Zhejiang Normal University, Jinhua, Zhejiang, 321004, P.R. China}}\\
		{\small{\it $^b$Xingzhi College, Zhejiang Normal University, Jinhua, Zhejiang, 321004, P.R. China}}
	}\maketitle
	\begin{abstract}
	Given a graph $G$ on $n$ vertices, its adjacency matrix and degree diagonal matrix are represented by $A(G)$ and $D(G)$, respectively. The $Q$-spectrum of $G$ consists of all the eigenvalues of its signless Laplacian matrix $Q(G)=A(G)+D(G)$ (including the multiplicities). A graph $G$ is known as being determined by its generalized $Q$-spectrum ($DGQS$ for short) if, for any graph $H$, $H$ and $G$ have the same $Q$-spectrum and so do their complements, then $H$ is isomorphic to $G$. In this paper, we present a method to construct $DGQS$ graphs. More specifically, let the matrix $W_{Q}(G)=\left [e,Qe,\dots ,Q^{n-1}e  \right ]$ ($e$ denotes the all-one column vector ) be the $Q$-walk matrix of $G$. It is shown that $G\circ P_{k}$ ($k=2,3$) is $DGQS$ if and only if $G$ is $DGQS$ for some specific graphs. This also provides a way to construct $DGQS$ graphs with more vertices by using $DGQS$ graphs with fewer vertices. At the same time, we also prove that $G\circ P_{2}$  is still $DGQS$ under specific circumstances. In particular, on the basis of the above results, we obtain an infinite sequences of $DGQS$ graphs $G\circ P_{k}^{t}$ ($k=2,3;t\ge 1$) for some specific $DGQS$ graph $G$.
	
	
	\end{abstract}

    \section{Introduction}
   All graphs in this paper are simple, that is, finite undirected graphs without multiple edges and loops. For a graph $G$ with $(0,1)$-adjacency matrix $A(G)$ and degree diagonal matrix $D(G)$, the matrix $Q(G)=D(G)+A(G)$ is called the signless Laplacian matrix of $G$ ($Q$-matrix for short). It is well-known that the $Q$-eigenvalues are the eigenvalues of $Q(G)$.
    
    Given a graph $G$ on $n$ vertices and a graph $H$ with a root vertex $u$, then the rooted product graph of $G$ and $H$ is obtained by copying one graph $G$ and copying $n$ graphs $H$, by gluing the $i$-th vertex of $G$ and the rooted vertex $u$ in the $i$-th copy of $H$ for $1\le i\le n$ (see Godsil and McKay \cite{Godsil1978}). Let $P_{k}$ be the rooted path of order $k$ and the root vertex be an endpoint, then the rooted product $G\circ P_{k}$ of graph $G$ and $P_{k}$ is described in Figure 1.
    
    Two graphs are cospectral if they have the same spectrum. A graph $G$ is known as being determined by its spectrum ($DS$ for short) if, for any graph $H$, $H$ and $G$ have the same spectrum, then $H$ is isomorphic to $G$. The generalized spectrum of $G$ is the spectrum of $G$ and the spectrum of its complement graph $\overline{G} $, where $\overline{G} $ is the complement of $G$ satisfying $A(\overline{G})=J-I-A(G)$.
    
    ``Which graphs are $DS$?" is an interesting question in graph theory. The problem goes back sixty years and is rooted in chemical theory (see \cite{Hs.H1956}). Another similar question was also raised by Kay in \cite{M.Kac1966}: ``Can one hear the shape of a drum?" This problem simulates the sound of the drum as the eigenvalue of the graph and the shape of the drum as the graph. However, it turns out that determining whether a graph is $DS$ is a difficult problem.
    
    It is well-known that $DS$ graphs are very few and have special structural properties. Wang and Xu in \cite{Wang20061,Wang20062} considered the above problem relative to the generalized spectrum determined in the controllable graph, where a graph $G$ is controllable if and only if the walk-matrix $W(G)=\left [e,Ae,\dots,A^{n-1}e\right ]$ is non-singular. A graph $G$ is known as being determined by its generalized spectrum ($DGS$ for short) if, for any graph $H$, $H$ and $G$ have the same spectrum and so do their complements, then $H$ is isomorphic to $G$. In $2006$, they also proposed a simple arithmetic condition to prove that a class of graphs is $DGS$ in which $G$ satisfies $\frac{{\rm det}(W(G))}{2^{\left \lfloor \frac{n}{2}  \right \rfloor } } $ is odd and square-free. Recently, some families of $DGS$ graphs were constructed by using the rooted product operation in \cite{Mao2022,Mao2015}. More specifically, it was proved that the rooted product $G\circ P_{k}$ ($k=2,3,4$) is $DGS$ when the graph $G$ satisfies ${\rm det}W(G)=\pm 2^{\frac{n}{2} } $ ($n$ is even) and $a_{0}=\pm 1$ (where $a_{0}$ is the constant term with respect to the characteristic polynomial of $G$). In \cite{Qiu2019}, the authors also presented a similar simple arithmetic condition to prove that a class of graphs is determined by its generalized $Q$-spectrum ($DGQS$ for short) in which $G$ satisfies $\frac{{\rm det}(W_{Q} (G))}{2^{\left \lfloor \frac{3n-2}{2}  \right \rfloor } }$ is odd and square-free. In this paper, we follow this research direction and explore how to construct $DGQS$ graphs.\\
    
    The contributions of this article are as follows:
    
    \begin{enumerate}[(1)]
    	\item We first establish an important determinant relationship between $W_{Q}(G)$ and $W_{Q}(G\circ P_{k})$ ($k=2,3$) as follows
    	\begin{equation}
    		{\rm det} W_{Q}(G\circ P_{2} ) =\pm a_{0}({\rm det} W_{Q}(G)) ^{2} ,~{\rm det} W_{Q}(G\circ P_{3} ) =\pm a_{0}^{2} ({\rm det} W_{Q}(G)) ^{3} ,
    	\end{equation}
        where $a_{0}$ is the constant term with respect to the characteristic polynomial of the matrix $Q(G)$.
    	
    	\item Based on the results above, we prove that the rooted product graph  $G\circ P_{k}$ ($k=2,3$) is $DGQS$ when ${\rm det}W_{Q}(G)=\pm 2^{\frac{3n-2}{2} } $ ($n$ is even) and $a_{0}=\pm2$. In particular, we also prove that each graph in $G\circ P_{k}^{t} $ ($k=2,3$ and $t\ge 1$) is $DGQS$.
    	
    	\item We study a family of graph $G$ satisfying  $\frac{{\rm det}W_{Q}(G) }{2^{\left \lfloor \frac{3n-2}{2}  \right \rfloor } } $ is odd and square-free. We show that the rooted product graph $G\circ P_{2}$ is still $DGQS$ in some cases.	
    \end{enumerate}
    
    The paper is organized as follows. In Section 2, we present some preliminary results. In Section 3, we give a relationship between ${\rm det}W_{Q}(G)$ and ${\rm det}W_{Q}(G\circ P_{k} )$ ($k=2,3$). Section 4 constructs some families of infinite $DGQS$ graphs. In Section 5, we show that $G\circ P_{2}$ is $DGQS$ under certain circumstances.
    
    \section{Preliminaries}
    In this section, we present some basic knowledges that will be used in the following sections. First,  we give two matrices $A$ and $B$ as below:
    
    $$A=\begin{bmatrix}
    	2&  1&  0&\cdots   &0  &0  &0 \\
    	1&  2&  1&  \cdots&  0&0  & 0\\
    	0&1  &2  &\cdots  &0  &0  &0 \\
    	\vdots & \vdots  & \vdots  & \ddots  & \vdots  &\vdots   &\vdots  \\
    	0&0  &0  &\cdots  &1  & 2 & 1\\
    	0&  0&  0& \cdots &0  &1  &1
    \end{bmatrix}_{(k-1)\times (k-1)},$$
    
    $$B=\begin{bmatrix}
    	1&1&  0&  0&\cdots   &0  &0  &0 \\
    	1&2&  1&  0&\cdots   &0  &0  &0 \\
    	0&1&  2&  1&  \cdots&  0&0  & 0\\
    	0&0&1  &2  &\cdots  &0  &0  &0 \\
    	\vdots & \vdots  & \vdots  & \vdots  & \ddots  &\vdots   &\vdots  &\vdots\\
    	0&0&0  &0  &\cdots  &1  & 2 & 1\\
    	0&0&  0&  0& \cdots &0  &1  &1
    \end{bmatrix}_{k\times k} .$$

    Suppose that the characteristic polynomials of matrix $A$ and matrix $B$ are $a_{k-1}(t) $ and $b_{k}(t)$, respectively. By a direct calculation, $a_{k-1}(t) $ and $b_{k}(t)$ satisfy the recursive relations of the following polynomials:
    \begin{equation}\label{equation1}
	a_{0}(t)=1,\quad a_{1}(t)=t-1 \quad {\rm and} \quad a_{k-1}(t)=(t-2)a_{k-2}(t)-a_{k-3}(t) \quad {\rm for} \quad k\ge3,
   \end{equation}
   \begin{equation}
	b_{0}(t)=1,\quad b_{1}(t)=t-1 \quad {\rm and} \quad b_{k}(t)=(t-1)a_{k-1}(t)-a_{k-2}(t) \quad {\rm for} \quad k\ge2.
   \end{equation}

    Next, we introduce the definition of rooted product graph. Let $G$ be a connected graph with vertex set $V(G)=\left \{ u_{1} ,u_{2},\dots ,u_{n} \right \}$ and $P_{k} $ is a path of length $k$ whose vertex set $V(P_{k})=\left \{ v_{1} ,v_{2},\dots ,v_{k} \right \}$. As shown in Figure 1, the rooted product graph of $G$ and $P_{k}$ is represented as $\widehat{G} _{k}=G\circ P_{k} $. It is obvious that $\widehat{G} _{k}$ has vertex set $V(\widehat{G}_{k}  )=\left \{ \left ( u_{i},v_{j}   \right ) \mid 1\le i\le n,1\le j\le k \right \}$ and edge set $E ( \widehat{G}_{k} ) =\left \{ \left ( u_{i},v_{1}   \right )\sim (u_{j},v_{1}  )\mid u_{i} u_{j}\in E(G)   \right \} \cup \left \{(u_{i},v_{s}  ) \sim (u_{i},v_{s+1}  )\mid 1\le i\le n, 1\le s\le k-1  \right \} $.
    \begin{figure}[h]
    	\centering
    	\includegraphics[width=10cm,height=6cm]{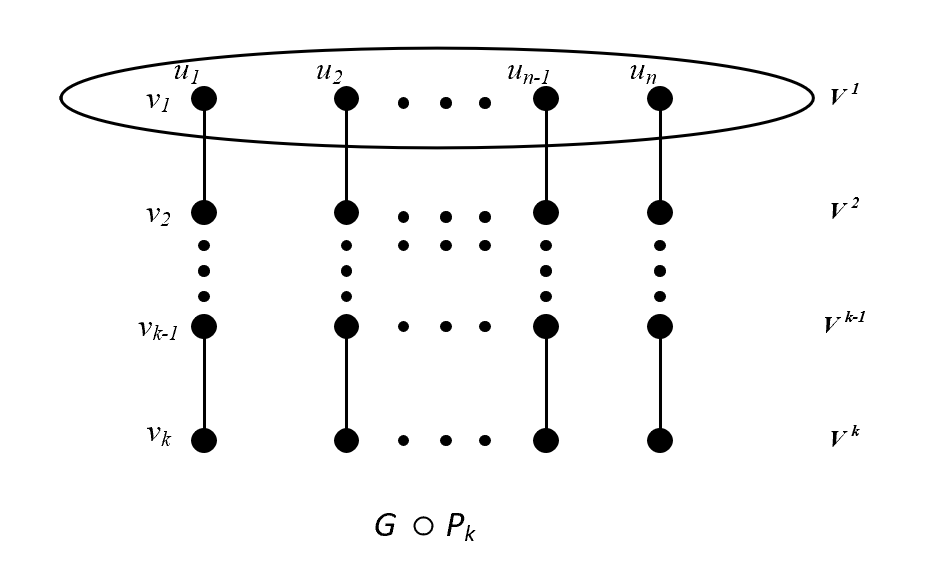}
    	\caption{The rooted product graph $\widehat{G}_{k}=G\circ P_{k}$.}
    	\label{Figure 1.}
    \end{figure}

    Let $V^{j}=\left \{\left ( u_{i},v_{j}   \right )\mid  1\le i\le n  \right \}$ for $1\le j\le k$, then $V^{1} ,V^{2} , \dots,V^{k} $ is a vertex partition of $V(\widehat{G}_{k})$. Following this partition, we can write the block matrix of the adjacency matrix of $\widehat{G}_{k}$ as follows:
    
    \begin{equation}
    	A(\widehat{G}_{k}  )=\begin{bmatrix}
    		A(G)& I_{n}  &0  &\cdots   & 0 & 0 & 0\\   
    		I_{n} &0  & I_{n} &\cdots   &0  & 0 & 0\\
    		0&  I_{n}&  0&  \cdots & 0 & 0 &0 \\
    		\vdots &\vdots   &\vdots   & \ddots  &\vdots   & \vdots  &\vdots  \\
    		0&  0&  0&  \cdots &  0& I_{n} &0 \\
    		0& 0 &  0&  \cdots &  I_{n}& 0 &I_{n} \\
    		0&  0&  0&  \cdots &  0& I_{n} &0
    	\end{bmatrix}.
    \end{equation}
    
    We know that the degree diagonal matrix of $\widehat{G}_{k}$ is $D(\widehat{G}_{k})={\rm diag} (D(G)+I_{n},2I_{n},\dots ,2I_{n} ,I_{n})$. Then the $Q$-matrix with respect to $\widehat{G}_{k}$ can be written as follows:
    
    \begin{equation}\label{eq5555}
    	Q(\widehat{G}_{k}  )=\begin{bmatrix}
    		Q(G)+I_{n}& I_{n}  &0  &\cdots   & 0 & 0 & 0\\   
    		I_{n} &2I_{n}  & I_{n} &\cdots   &0  & 0 & 0\\
    		0&  I_{n}&  2I_{n}&  \cdots & 0 & 0 &0 \\
    		\vdots &\vdots   &\vdots   & \ddots  &\vdots   & \vdots  &\vdots  \\
    		0&  0&  0&  \cdots &  2I_{n}& I_{n} &0 \\
    		0& 0 &  0&  \cdots &  I_{n}& 2I_{n} &I_{n} \\
    		0&  0&  0&  \cdots &  0& I_{n} &I_{n}
    	\end{bmatrix}.
    \end{equation}

    In \cite{Mao2015}, the authors gave an interesting relationship between ${\rm det}W(G)$ and ${\rm det}W(G\circ P_{2})$ as described below.
    
    \paragraph{Theorem 2.1}(\cite{Mao2015}). Let the graph $G\circ P_{2}$ be obtained by taking the rooted product of a graph $G$ and the graph $P_{2}$. Then ${\rm det} W(G\circ P_{2} )=\pm a_{0}({\rm det}  (W(G)))^{2} $, where $a_{0}$ is the constant term with respect to the characteristic polynomial of $G$.\\
    
    In \cite{Mao2022}, the authors also gave the relation between ${\rm det}W(G)$ and ${\rm det}W(G\circ P_{k})$ $(k=3,4)$ as described below.
    
    \paragraph{Theorem 2.2} (\cite{Mao2022}). Let the graph $G\circ P_{3}$ (resp. $G\circ P_{4}$) be obtained by taking the rooted product of a graph $G$ and the graph $P_{3}$ (resp. $P_{4}$). Then ${\rm det} W(G\circ P_{3} )=\pm a_{0}({\rm det}(W(G)))^{3} $ (resp. ${\rm det}W(G\circ P_{4})=\pm a_{0}^{2} ({\rm det}(W(G)) )^{4}$, where $a_{0}$ is the constant term with respect to the characteristic polynomial of $G$.\\
    
    In \cite{Mao2022},the authors proposed a conjecture, which is proved in \cite{Wang2023} as follows.
    
    \paragraph{Theorem 2.3}(\cite{Wang2023}). Let the graph $G\circ P_{k}$ be obtained by taking the rooted product of a graph $G$ and the graph $P_{k}$. Then ${\rm det} W(G\circ P_{k} )=\pm a_{0}^{\left \lfloor \frac{k}{2}  \right \rfloor }({\rm det}(W(G)))^{k} $, where $a_{0}$ is the constant term with respect to the characteristic polynomial of $G$.\\
    
    In this paper, we replace $W(G)$ with $W_{Q}(G)$, then a similar relationship is obtained as follows.
    
    \paragraph{Theorem 2.4.} Let the graph $G\circ P_{2}$ (resp. $G\circ P_{3}$) be obtained by taking the rooted product of a graph $G$ and the graph $P_{2}$ (resp. $P_{3}$). Then ${\rm det} W_{Q}(G\circ P_{2} )=\pm a_{0}({\rm det}  (W_{Q}(G)))^{2} $ (resp. ${\rm det}W_{Q}(G\circ P_{3})=\pm a_{0}^{2} ({\rm det}(W_{Q}(G)))^{3}$, where $a_{0}$ is the constant term with respect to the characteristic polynomial of the matrix $Q(G)$.\\
    
    Before proving Theorem 2.4, we need the following lemma.
    
    \paragraph{Lemma 2.5}\label{le2.1} (\cite{Lou20172}). Let $\widehat{G}_{k}=G\circ P_{k}$ with the $Q$-matrix in \eqref{eq5555}. If $t$ is a $Q$-eigenvalue of $\widehat{G}_{k}$ with the corresponding eigenvector
    \begin{equation*}
    	\zeta =\left [ x_{1}^{T} , x_{2}^{T},\dots ,x_{k}^{T}\right ] ^{T},
    \end{equation*}
    where the vertices of $V^{s}$ correspond to elements of $x_{s}$ for $1\le s \le k $, then $\lambda =\frac{a_{k}(t) }{a_{k-1}(t) }  +1$ is a $Q$-eigenvalue of $G$ with the corresponding eigenvector $x_{1}$, and $x_{s}=\frac{a_{k-s}(t) }{a_{k-1}(t) }x_{1} $ ($a_{k-1}(t)\ne 0$) for $2\le s\le k$, where $a_{s}(t)$ is defined in \eqref{equation1}.
    
    \paragraph{Remark 2.6.} Due to Lemma 2.5, we know that the $Q$-eigenvalue $\lambda$ of $G$ satisfies $\lambda =\frac{a_{k}(t) }{a_{k-1}(t) }  +1$. Using the recursive relation of the polynomial $a_{k-1}(t)$, we can get that $\lambda =\frac{a_{k}(t)+a_{k-1}(t)  }{a_{k-1}(t) } =\frac{(t-2)a_{k-1}(t)-a_{k-2}(t)+a_{k-1}(t)   }{a_{k-1}(t) } =\frac{b_{k}(t) }{a_{k-1}(t) }$. Then the $Q$-eigenvalue $t$ of  $\widehat{G}_{k}$ satisfies the following polynomial: $$\varphi(t)=b_{k}(t) - \lambda a_{k-1}(t). $$
    Since $\varphi(t)$ is a polynomial of degree $k$, we can conclude that one $Q$-eigenvalue of $G$ corresponds exactly to $k$ $Q$-eigenvalues of $\widehat{G}_{k}$.\\
    
    The following theorem is an immediate conclusion from Lemma 2.5 and Remark 2.6.
    
    \paragraph{Theorem 2.7.} Let $\lambda _{i} $ be the eigenvalue of the graph $G$ with respect to the normalized eigenvector $x_{i}$ ($i=1,2,\dots,n$). Let $t_{i}^{(j)} $ ($i=1,2,\dots,n;j=1,2,\dots,k$  ) be the $Q$-eigenvalues of the rooted product graph $\widehat{G}_{k}$, and the eigenvector $\zeta _{i}^{j}$ corresponding to the $Q$-eigenvalue  $t_{i}^{(j)} $ can be written $\zeta _{i}^{j} =\left [ x_{i}^{T},\frac{a_{k-2}(t_{i}^{(j)} ) }{a_{k-1}(t_{i}^{(j)} ) }x_{i}^{T} ,\dots , \frac{a_{0}(t_{i}^{(j)} ) }{a_{k-1}(t_{i}^{(j)} ) }x_{i}^{T}    \right ] ^{T} $. In addition, $t_{i}^{(j)}$ satisfies the following equation:
    $$\varphi(t_{i}^{(j)} )=b_{k}(t_{i}^{(j)} ) - \lambda_{i} a_{k-1}(t_{i}^{(j)} ),$$
    where $b_{k}(t)$ and $a_{k-1}(t)$ are characteristic polynomials of $B$ and $A$, respectively.
    
    \paragraph{Lemma 2.8}(\cite{Vinberg2003}). Any symmetric polynomial can be uniquely expressed as a polynomial in elementary symmetric polynomials.
    
    \paragraph{Definition 2.9.} Given any two polynomials $f(x)=a_{n}x^{n}+a_{n-1}x^{n-1}+\cdots +a_{1}x+a_{0}\in \mathbb{R} \left [ x \right ] $ and $g(x)=b_{m}x^{m}+b_{m-1}x^{m-1}+\cdots +b_{1}x+b_{0}\in \mathbb{R} \left [ x \right ] $, the resultant of $f(x)$ and $g(x)$ is expressed as Res($f(x),g(x)$), which is defined to be
    $$a_{n}^{m}b_{m}^{n} \prod_{1\le i\le n,1\le j\le m}(\alpha _{i}-\beta _{j}  ),$$
    where $\alpha _{i}$$^{,}$s $\in \mathbb{C} $ are the root of $f(x)$ and $\beta_{i}$$^{,}$s $\in \mathbb{C} $ are the root of $g(x)$.

    \section{The determinant of $W_{Q}(\widehat{G}_{k})$}
    In this section, We shall show an interesting relationship between the determinants of $W_{Q}(G)$ and the determinants of $W_{Q}(\widehat{G}_{k})$. First of all, we will prove that the determinant of $W_{Q}(G)$ is determined by the eigenvalues and eigenvectors of the matrix $Q(G)$.
    
    \paragraph{Lemma 3.1.}\label{le3.3} Let $\lambda _{i} $ be the $Q$-eigenvalue of $G$ with the corresponding normalized eigenvector $x_{i} $ for $i=1,2,\dots,n$. Then 
    $${\rm det}W_{Q}(G) =\pm \prod_{1\le i< j\le n} (\lambda _{i}-\lambda _{j}  )\prod_{1\le i\le n} (e_{n}^{T} x_{i} ).$$
    \begin{proof}
    	On the grounds of $Q(G)x_{i}=\lambda _{i}x_{i}$ for $i=1,2,\dots,n$, then $Q^{k}(G) x_{i}=\lambda _{i}^{k} x_{i}$ for $i=1,2,\dots,n$. Since
    	\begin{equation*}
    		\begin{aligned}
    			\begin{bmatrix}
    				e^{T} \\
    				e^{T}Q(G) \\
    				\vdots \\
    				e^{T}Q^{n-1}(G ) 
    			\end{bmatrix}\left [ x_{1},x_{2,\dots },x_{n}    \right ] &=\begin{bmatrix}
    				e^{T}x_{1}   &e^{T}x_{2}  & \cdots  &e^{T}x_{n}  \\
    				\lambda _{1} e^{T}x_{1} & \lambda _{2} e^{T}x_{2} & \cdots  & \lambda _{n} e^{T}x_{n}\\
    				\vdots & \vdots  & \ddots  &\vdots  \\
    				\lambda _{1}^{n-1} e^{T}x_{1} &\lambda _{2}^{n-1} e^{T}x_{2}  & \cdots  &\lambda _{n}^{n-1} e^{T}x_{n}
    			\end{bmatrix}\\
    			&=\begin{bmatrix}
    				1&  1&  \cdots & 1\\
    				\lambda _{1} & \lambda _{2} & \cdots  & \lambda _{n}\\
    				\vdots & \vdots  & \ddots  &\vdots  \\
    				\lambda _{1}^{n-1} & \lambda _{2}^{n-1} & \cdots  &\lambda _{n}^{n-1}
    			\end{bmatrix}\begin{bmatrix}
    				e^{T}x_{1}  &  &  & \\
    				&  e^{T}x_{2}&  & \\
    				&  &  \ddots & \\
    				&  &  &e^{T}x_{n}
    			\end{bmatrix}.
    		\end{aligned}
    	\end{equation*}
Then, by taking the determinants of both sides of the above equation and ${\rm det}(\left [ x_{1},x_{2},\dots,x_{n}\right ] )=\pm 1$, we obtain the required result.
    \end{proof}
    
For the sake of convenience, let $f_{k}(t)=a_{0}(t)+a_{1}(t)+\cdots +a_{k-1}(t)$, which will be useful in the following proof.
    
    \paragraph{Theorem 3.2.} Let $\lambda _{i} $ be the $Q$-eigenvalues of the graph $G$ with respect to the normalized eigenvectors $x_{i}$ ($i=1,2,\dots ,n$). Let $t_{i}^{(j)}$ be the eigenvalues of the rooted product graph $G\circ P_{k}$ with respect to the eigenvectors $\zeta _{i}^{j}$ ($i=1,2,\dots,n;~j=1,2,\dots ,k$). Then
    \begin{equation*}
    	{\rm det}(W_{Q}(G\circ P_{k} ))=\pm({\rm det}(W_{Q}(G))^{k}\prod_{1\le i\le n} \prod_{1\le k_{2}\le k }a_{k-1}(t_{i}^{(k_{2} )} ) f_{k}(t_{i}^{(k_{2}) } ).
    \end{equation*}
    \begin{proof}
        Suppose that the signless Laplacian matrix of rooted product graph $G\circ P_{k}$ is $\widetilde{Q} =Q(G\circ P_{k} ) $. Then we can obtain, for $j=1,2,\dots ,k,$
        \begin{equation}\label{eqeq5}
        	\begin{bmatrix}
        		e_{kn}^{T} \\
        		e_{kn}^{T}\widetilde{Q} \\
        		\vdots \\
        		e_{kn}^{T}\widetilde{Q}^{kn-1}  
        	\end{bmatrix}\left [ E^{(1)},E^{(2)},\dots ,E^{(k)}  \right ] =\left [ M^{(1)},M^{(2)},\dots ,M^{(k)}  \right ]D_{kn},
        \end{equation}
        where 
        $$E^{(j)}= \left [\zeta _{1}^{j},  \zeta _{2}^{j},\dots ,\zeta _{n}^{j} \right ] ,$$
        $$ M^{(j)}=\begin{bmatrix}
        	1 & 1 &  \cdots &1 \\
        	t_{1}^{(j)} &  t_{2}^{(j)} & \cdots  & t_{n}^{(j)} \\
        	(t_{1}^{(j)})^{2} &(t_{2}^{(j)})^{2}  & \cdots  &(t_{n}^{(j)})^{2} \\
        	\vdots & \vdots  & \ddots   & \vdots \\
        	(t_{1}^{(j)})^{kn-1}& (t_{2}^{(j)})^{kn-1} & \cdots  &(t_{n}^{(j)})^{kn-1}
        \end{bmatrix},$$
        $$D_{kn}={\rm diag}\left( \frac{f_{k}(t_{1}^{(1)} ) }{a_{k-1}(t_{1}^{(1)} ) }e_{n}^{T}x_{1} , \frac{f_{k}(t_{2}^{(1)} ) }{a_{k-1}(t_{2}^{(1)} ) }e_{n}^{T}x_{2} ,\dots ,\frac{f_{k}(t_{n}^{(1)} ) }{a_{k-1}(t_{n}^{(1)} ) }e_{n}^{T}x_{n} ,\dots,\right.$$ $$\quad \left.\frac{f_{k}(t_{1}^{(k)} ) }{a_{k-1}(t_{1}^{(k)}) } e_{n}^{T}x_{1}  ,\frac{f_{k}(t_{2}^{(k)} ) }{a_{k-1}(t_{2}^{(k)}) } e_{n}^{T}x_{2} ,\dots ,\frac{f_{k}(t_{n}^{(k)} ) }{a_{k-1}(t_{n}^{(k)}) } e_{n}^{T}x_{n} \right).$$
        Taking the determinant of both sides of the formula $(6)$, we obtain
        $${\rm det}(W_{Q}(G\circ P_{k} )){\rm det} \left [ E^{(1)}, E^{(2)},\dots ,E^{(k)} \right ]={\rm det}(\left [ M^{(1)}, M^{(2)},\dots ,M^{(k)} \right ] )\prod_{1\le i\le n}(e_{n}^{T}x_{i}  )^{k}\prod_{1\le i\le n} \prod_{1\le k_{2} \le k}\frac{f_{k}(t_{i}^{(k_{2} )} ) }{a_{k-1}(t_{i}^{(k_{2}) } ) }.$$
        
        Remark that, in the following calculation, we only consider the absolute value of the determinant. First of all, we have
        \begin{equation*}
        	\begin{aligned}
        		 &{\rm det}\left [ E^{(1)}, E^{(2)},\dots ,E^{(k)}   \right ] \\
        		 &={\rm det}\left [ \zeta _{1}^{1}, \zeta _{2}^{1},\dots ,\zeta _{n}^{1},\dots ,\zeta _{1}^{k},\zeta _{2}^{k},\dots , \zeta _{n}^{k}\right ] \\ 
        		&={\rm det}\begin{bmatrix}
        			x_{1} & x_{2}  & \cdots  & x_{n}  & \cdots  & x_{1}  & x_{2}  &\cdots   &x_{n}  \\
        			\frac{a_{k-2}(t_{1}^{(1)} ) }{a_{k-1}(t_{1}^{(1)} ) }x_{1}  &\frac{a_{k-2}(t_{2}^{(1)} ) }{a_{k-1}(t_{2}^{(1)} ) }x_{2}  & \cdots  &\frac{a_{k-2}(t_{n}^{(1)} ) }{a_{k-1}(t_{n}^{(1)} ) }x_{2}  &\cdots   & \frac{a_{k-2}(t_{1}^{(k)} ) }{a_{k-1}(t_{1}^{(k)} ) }x_{1}  &  \frac{a_{k-2}(t_{2}^{(k)} ) }{a_{k-1}(t_{2}^{(k)} ) }x_{2} & \cdots  &  \frac{a_{k-2}(t_{n}^{(k)} ) }{a_{k-1}(t_{n}^{(k)} ) }x_{n}\\
        			\vdots & \vdots  &\ddots   &\vdots   & \ddots   & \vdots  & \vdots  &\ddots    & \vdots \\
        			\frac{1 }{a_{k-1}(t_{1}^{(1)} ) }x_{1}& \frac{1 }{a_{k-1}(t_{2}^{(1)} ) }x_{x} &\cdots   &\frac{1 }{a_{k-1}(t_{n}^{(1)} ) }x_{n}  & \cdots  & \frac{1 }{a_{k-1}(t_{1}^{(k)} ) }x_{1} &\frac{1 }{a_{k-1}(t_{2}^{(k)} ) }x_{2}  &\cdots   &\frac{1 }{a_{k-1}(t_{n}^{(k)} ) }x_{n}
        		\end{bmatrix} \\
        		&={\rm det}\begin{bmatrix}
        			x_{1}  & x_{2}  & \cdots  & x_{n}  &  &  &  &  & \\
        			&  &  &  & \ddots  &  &  &  & \\
        			&  &  &  &  &  x_{1}  & x_{2}  & \cdots  & x_{n}
        		\end{bmatrix}{\rm det}L,
        	\end{aligned}
        \end{equation*}
       where
       
        $$L=\begin{bmatrix}
        	1& 0& \cdots  & 0 & \cdots  & 1 & 0 &\cdots   &0 \\
        	0&  1& \cdots  & 0 &\cdots   &0  &1  &\cdots   &0 \\
        	\vdots & \vdots  &\ddots  &0  &\cdots   &\vdots   &\vdots   &\ddots   &\vdots  \\
        	0& 0 & \cdots  &1  &\cdots  & 0 & 0 &\cdots   & 1\\
        	\vdots & \vdots  & \vdots  &\vdots   & \ddots  &\vdots   &\vdots   &\vdots   &\vdots  \\
        	\frac{1}{a_{k-1}(t_{1}^{(1)} )} &0  &\cdots   & 0 &\cdots   &\frac{1}{a_{k-1}(t_{1}^{(k)} )}  & 0 & \cdots  &0 \\
        	0&  \frac{1}{a_{k-1}(t_{2}^{(1)} )}&  \cdots &  0& \cdots  & 0 & \frac{1}{a_{k-1}(t_{2}^{(k)} )}  & \cdots &0 \\
        	\vdots & \vdots  & \ddots &\vdots   & \vdots  &\vdots    &\vdots   & \ddots  &\vdots  \\
        	0& 0 & \cdots  &  \frac{1}{a_{n}(t_{1}^{(1)} )}&\cdots   &0  & 0 & \cdots  &\frac{1}{a_{k-1}(t_{n}^{(k)} )}
        \end{bmatrix}.$$
        We first extract the common factor of the matrix $L$ and then do the elementary row transformation of the matrix $L$,
        \begin{equation*}
        \begin{aligned}
        		{\rm det}\left [ E^{(1)}, E^{(2)},\dots ,E^{(k)}   \right ]&=\prod_{1\le i\le n}(\frac{1}{\prod_{1\le k_{2} \le k}a_{k-1}(t_{i}^{(k_{2}) } ) }) {\rm det}\begin{bmatrix}
        			1 & 1 & \cdots  &1 \\
        			t_{i}^{(1)} & t_{i}^{(2)} & \cdots  &t_{i}^{(k)} \\
        			\vdots &  \vdots & \ddots  & \vdots \\
        			(t_{i}^{(1)})^{k-1} & (t_{i}^{(2)})^{k-1}  & \cdots  &(t_{i}^{(k)})^{k-1} 
        		\end{bmatrix}\\
        		&=\frac{\prod_{1\le i\le n}\prod_{1\le k_{1}< k_{2}\le k  }(t_{i}^{(k_{2}) }-t_{i}^{(k_{1}) }  ) }{\prod_{1\le i\le n}\prod_{1\le k_{2}\le k  }(a_{k-1}(t_{i}^{(k_{2}) } ) )}.
        \end{aligned}   	
        \end{equation*}
        
        Next, we compute the Vandermonde determinant ${\rm det}\left [ E^{(1)}, E^{(2)},\dots ,E^{(k)} \right ].$ Assume that ${\rm det} \left [ E^{(1)}, E^{(2)},\dots ,E^{(k)} \right ] =V_{1}V_{2},$ where 
        \begin{equation*}
        \begin{aligned}
        	    V_{1}&=\prod_{1\le i\le n} \prod_{1\le k_{1}<k_{2}\le k}(t_{i}^{(k_{2} )}-t_{i}^{(k_{1}) }  ),\\
        		V_{2}&=\prod_{1\le i< j\le n} \prod_{1\le k_{2}\le k }\prod_{1\le k_{1}\le k }(t_{j}^{(k_{2}) } -t_{i}^{(k_{1}) }  ) \\
        		&=\prod_{1\le i< j\le n} \prod_{1\le k_{2}\le k }(b_{k}(t_{j}^{(k_{2}) } ) - \lambda_{i} a_{k-1}(t_{j}^{(k_{2}) } ))\\
        		&=\prod_{1\le i< j\le n} \prod_{1\le k_{2}\le k }(\lambda _{j}a_{k-1}(t_{j}^{(k_{2}) } )  - \lambda_{i}  a_{k-1}(t_{j}^{(k_{2}) } ))\\
        		&=\prod_{1\le i< j\le n} (\lambda_{j}-  \lambda_{i}  )^{k}\prod_{1\le k_{2}\le k }a_{k-1} (t_{j}^{(k_{2}) } ).\\
        \end{aligned}
        \end{equation*}
        Thus, we can get that
        \begin{equation*}
        \begin{aligned}
        		{\rm det}(W_{Q}(G\circ P_{k} ))&=\prod_{1\le i\le n}(e_{n}^{T}x_{i}  )^{k}\prod_{1\le i< j\le n}(\lambda _{j}-\lambda _{i}  )^{k}\prod_{1\le j\le n}\prod_{1\le k_{2}\le k }a_{k-1}(t_{j}^{(k_{2} )} ) \prod_{1\le i\le n} \prod_{1\le k_{2}\le k }f_{k}(t_{j}^{(k_{2}) } )\\
        		&=({\rm det}(W_{Q}(G))^{k}\prod_{1\le j\le n} \prod_{1\le k_{2}\le k }a_{k-1} (t_{j}^{(k_{2} )})\prod_{1\le i\le n}\prod_{1\le k_{2} \le k} f_{k}(t_{i}^{(k_{2} )} )\\
        		&=({\rm det}(W_{Q}(G))^{k}\prod_{1\le i\le n} \prod_{1\le k_{2}\le k }a_{k-1}(t_{i}^{(k_{2} )} ) f_{k}(t_{i}^{(k_{2}) } ) .\\
       	\end{aligned}
        \end{equation*}
    Since we only consider the absolute values in the process of the proof, then the conclusion is clearly valid. This completes the proof.
    \end{proof}
    
    In what follows, we shall give the proof of Theorem 2.4:
    
    \paragraph{Proof of Theorem 2.4.} 
    \begin{proof}
    	We only prove the case of $G\circ P_{3}$ in Theorem 2.4, the case of $G\circ P_{2}$ can be proved in the same way.
    	From Theorem 3.2, we know $${\rm det}(W_{Q}(G\circ P_{k})))=\pm({\rm det}(W_{Q}(G))^{k}\prod_{1\le i\le n} \prod_{1\le k_{2}\le k }a_{k-1}(t_{i}^{(k_{2} )} ) f_{k}(t_{i}^{(k_{2}) } ).$$
    	When $k=3$, we have
    	$${\rm det}(W_{Q}(G\circ P_{3})))=\pm ({\rm det}(W_{Q}(G))^{3}\prod_{1\le i\le n} \prod_{1\le k_{2}\le 3 }a_{2}(t_{i}^{(k_{2} )} ) f_{3}(t_{i}^{(k_{2}) } ).$$
    	Notice that the roots of quadratic equation $a_{2}(t)$ are $\frac{3+\sqrt[]{5}}{2}$, $\frac{3-\sqrt[]{5}}{2} $, and the roots of quadratic equation $f_{3}(t)$ are $1$, $1$.
    	Then it can be concluded that $$\prod_{1\le i\le n} \prod_{1\le k_{2}\le 3 }a_{2}(t_{i}^{(k_{2} )} ) f_{3}(t_{i}^{(k_{2}) } )=\prod_{1\le i\le n} \prod_{1\le k_{2}\le 3 }(t_{i}^{(k_{2})}-\frac{3+\sqrt[]{5}}{2} )(t_{i}^{(k_{2})}-\frac{3-\sqrt[]{5}}{2} )(t_{i}^{(k_{2})}-1)(t_{i}^{(k_{2})}-1).$$
    	Now we shall compute $\prod_{1\le i\le n} \prod_{1\le k_{2}\le 3 }(t_{i}^{(k_{2})}-\frac{3+\sqrt[]{5}}{2} )$, $\prod_{1\le i\le n} \prod_{1\le k_{2}\le 3 }(t_{i}^{(k_{2})}-\frac{3-\sqrt[]{5}}{2} )$ and  $\prod_{1\le i\le n} \prod_{1\le k_{2}\le 3 }(t_{i}^{(k_{2})}-1)$ as below:

    	\begin{equation}
    	\begin{aligned}
    			 \prod_{1\le i\le n} \prod_{1\le k_{2}\le 3 }(t_{i}^{(k_{2})}-\frac{3+\sqrt[]{5}}{2})&=\prod_{1\le i\le n}(t_{i}^{(1)}-\frac{3+\sqrt[]{5}}{2})(t_{i}^{(2)}-\frac{3+\sqrt[]{5}}{2})(t_{i}^{(3)}-\frac{3+\sqrt[]{5}}{2})\\
    				&=t_{i}^{(1)}t_{i}^{(2)} t_{i}^{(3)}-\frac{3+\sqrt[]{5} }{2}(t_{i}^{(1)}t_{i}^{(2)}+ t_{i}^{(1)}t_{i}^{(3)}+t_{i}^{(2)} t_{i}^{(3)})\\
    				&\quad+\frac{7+3\sqrt[]{5} }{2}(t_{i}^{(1)}+t_{i}^{(2)}+ t_{i}^{(3)})-9-4\sqrt[]{5},
    	\end{aligned}
    	\end{equation}

        %\begin{small}
    	\begin{equation}
    	\begin{aligned}
    			\prod_{1\le i\le n} \prod_{1\le k_{2}\le 3 }(t_{i}^{(k_{2})}-\frac{3-\sqrt[]{5}}{2})&=\prod_{1\le i\le n}(t_{i}^{(1)}-\frac{3-\sqrt[]{5}}{2})(t_{i}^{(2)}-\frac{3-\sqrt[]{5}}{2})(t_{i}^{(3)}-\frac{3-\sqrt[]{5}}{2})\\
    			&=t_{i}^{(1)}t_{i}^{(2)} t_{i}^{(3)}-\frac{3-\sqrt[]{5} }{2}(t_{i}^{(1)}t_{i}^{(2)}+ t_{i}^{(1)}t_{i}^{(3)}+t_{i}^{(2)} t_{i}^{(3)})\\
    			&\quad+\frac{7-3\sqrt[]{5} }{2}(t_{i}^{(1)}+t_{i}^{(2)}+ t_{i}^{(3)})-9+4\sqrt[]{5},
    	\end{aligned}
    	\end{equation}
        %\end{small} 

    	\begin{equation}
    	\begin{aligned}
    			\prod_{1\le i\le n} \prod_{1\le k_{2}\le 3 }(t_{i}^{(k_{2})}-1)&=\prod_{1\le i\le n}(t_{i}^{(1)}-1)(t_{i}^{(2)}-1)(t_{i}^{(3)}-1)\\
    			&=t_{i}^{(1)}t_{i}^{(2)} t_{i}^{(3)}-(t_{i}^{(1)}t_{i}^{(2)}+ t_{i}^{(1)}t_{i}^{(3)}+t_{i}^{(2)} t_{i}^{(3)})+(t_{i}^{(1)}+t_{i}^{(2)}+ t_{i}^{(3)})-1.
    	\end{aligned}
        \end{equation}
Since $t$ is the root of the equation $$\varphi(t_{i}^{(j)} )=b_{3}(t_{i}^{(j)} ) - \lambda_{i} a_{2}(t_{i}^{(j)} )=(t_{i}^{(j)})^{3}-(4+\lambda _{i} )(t_{i}^{(j)})^{2} +(3+3\lambda _{i} )t_{i}^{(j)}-\lambda _{i},$$ 
then $$(t_{i}^{(1)}+t_{i}^{(2)}+ t_{i}^{(3)})=4+\lambda _{i},$$ 
    $$(t_{i}^{(1)}t_{i}^{(2)}+ t_{i}^{(1)}t_{i}^{(3)}+t_{i}^{(2)} t_{i}^{(3)})=3+3\lambda _{i},$$ $$t_{i}^{(1)}t_{i}^{(2)} t_{i}^{(3)}=\lambda _{i}.$$ By substituting the above equation into (7), (8) and (9), we get 
    $$\prod_{1\le i\le n} \prod_{1\le k_{2}\le 3 }(t_{i}^{(k_{2})}-\frac{3+\sqrt[]{5}}{2})=\frac{1+\sqrt[]{5} }{2},$$ 
    $$\prod_{1\le i\le n} \prod_{1\le k_{2}\le 3 }(t_{i}^{(k_{2})}-\frac{3-\sqrt[]{5}}{2})=\frac{1-\sqrt[]{5} }{2},$$ 
    $$\prod_{1\le i\le n} \prod_{1\le k_{2}\le 3 }(t_{i}^{(k_{2})}-1)=-\lambda _{i}.$$
    Hence, we get $${\rm det}(W_{Q}(G\circ P_{3} ))=\pm\prod_{1\le i\le n}(\lambda _{i})^{2}({\rm det}(W_{Q}(G))^{3}=\pm a_{0}^{2} ({\rm det}(W_{Q}(G)))^{3}.$$ This completes the proof.
    \end{proof}
    We conclude this section by asking a question similar to the conjecture in \cite{Mao2022}.
    
    \paragraph{Question 3.3.} Let $G$ be a graph with $a_{0}$ is the constant term with respect to the characteristic polynomial of the matrix $Q(G)$. Is it true that ${\rm det} W(G\circ P_{k} )=\pm a_{0}^{\left \lceil \frac{k}{2}  \right \rceil }({\rm det}(W(G)))^{k} $ whenever $k\ge 4$?

    \section{Constructig a infinite families of $DGQS$ graphs}
    In this section, we use the method of the rooted product operation to construct some families of infinite $DGQS$ graphs. For convenience, we need to define graph $G\circ P_{k} ^{t} $ ($k\ge 2,t\ge 1$). Graph $G\circ P_{k} ^{t} $ is the graph obtained by taking the rooted product of graph $G\circ P_{k} ^{t-1} $ and the graph $P_{k}$. For example, $G\circ P_{k}^{0}=G$, $G\circ P_{k}^{1}=G\circ P_{k}$, $G\circ P_{k}^{2}=(G\circ P_{k})\circ P_{k}$, $G\circ P_{k}^{3}=((G\circ P_{k})\circ P_{k})\circ P_{k}$,$\dots$.
    
    First, we recall some known results in \cite{Qiu2019}, the authors defined a family of graphs
    	$$\mathcal{F}_{Q,n}=\left \{ G:\frac{{\rm det}W_{Q}(G) }{2^{\left \lfloor \frac{3n-2}{2}  \right \rfloor } }  ~{\rm is~an~odd~and~square{\mbox -}free~integer}  \right \} $$
    and proved the following theorem.
    
    \paragraph{Theorem 4.1}(\cite{Qiu2019}). Every graph in $\mathcal{F}_{Q,n}$ is $DGQS$.\\
    
    In \cite{Mao2022}, the authors first showed the relationship between ${\rm det}W(G)$ and ${\rm det}W(G\circ P_{k})$ ($k=3,4$) and then proved the following two theorems.
    
    \paragraph{Theorem 4.2}(\cite{Mao2022}). Let $G$ be a graph with ${\rm det}W(G)=\pm 2^{\frac{n}{2} }$ ($n$ is even) and $a_{0}=\pm1$. Then $G\circ P_{k} $ ($k=3,4$) are $DGS$.
    
    \paragraph{Theorem 4.3}(\cite{Mao2022}). Let $G$ be a graph with ${\rm det}~W(G)=\pm 2^{\frac{n}{2} }$ ($n$ is even) and $a_{0}=\pm1$. Then every graph in the infinite sequence $G\circ P_{k}^{t} $ ($k=3,4;~t \ge 1$) is $DGS$.\\
    
    Next, we can get similar results to Theorems 4.2 and 4.3 by replacing the matrix $W(G)$ and $a_{0}=\pm1$ with $W_{Q}(G)$ and $a_{0}=\pm2$, respectively. Meanwhile, the proof of the following theorem is mainly based on Theorem 4.1.
    
    \paragraph{Theorem 4.4.} Let $G$ be a graph with ${\rm det}W_{Q}(G)=\pm 2^{\frac{3n-2}{2} } $ ($n$ is even) and $a_{0}=\pm2$. Then $G\circ P_{k} $ ($k=2,3$) are $DGQS$.
    \begin{proof}
    	Let $G$ be a graph of $n$ vertices. From the definition of the rooted product operation, we know that $G\circ P_{k} $ ($k=2,3$) have $kn$ vertices. Applying the result of Theorem 2.4, we can get ${\rm det} W_{Q}(G\circ P_{2} )=\pm a_{0}{\rm det}^{2}  (W_{Q}(G))=\pm 2^{3n-1}$ and ${\rm det}W_{Q}(G\circ P_{3})=\pm a_{0}^{2} {\rm det}^{3}(W_{Q}(G))=\pm 2^{\frac{9n-2}{2} } $. By a direct calculation, we get $\frac{{\rm det}W_{Q}(G\circ P_{2} ) }{2^{\frac{6n-2}{2} } } =\pm 1$ and $\frac{{\rm det}W_{Q}(G\circ P_{3} ) }{2^{\frac{9n-2}{2} } } =\pm 1$. Then we obtain  $G\circ P_{k} $ ($k=2,3$) are $DGQS$ from Theorem 4.1.
    \end{proof}	
    
    \paragraph{Lemma 4.5.}	Let $G$ be a graph with $a_{0}=\pm 2$, where $a_{0}$ is the constant term with respect to the characteristic polynomial of the matrix $Q(G)$. Also let $a_{0}^{(t)} $ be the constant term with respect to the characteristic polynomial of the matrix $Q(G\circ P_{k} ^{t})$. Then $a_{0}^{(t)}=\pm 2$.
    
    \begin{proof}
    	We will prove this lemma by induction on $t$. First of all, if $t= 1$, then
    	\begin{equation*}
    		\begin{aligned}
    			a_{0}^{(1)}&={\rm det}\begin{bmatrix}
    				Q(G)+I_{n}& I_{n}  &0  &\cdots   & 0 & 0 & 0\\   
    				I_{n} &2I_{n}  & I_{n} &\cdots   &0  & 0 & 0\\
    				0&  I_{n}&  2I_{n}&  \cdots & 0 & 0 &0 \\
    				\vdots &\vdots   &\vdots   & \ddots  &\vdots   & \vdots  &\vdots  \\
    				0&  0&  0&  \cdots &  2I_{n}& I_{n} &0 \\
    				0& 0 &  0&  \cdots &  I_{n}& 2I_{n} &I_{n} \\
    				0&  0&  0&  \cdots &  0& I_{n} &I_{n}
    			\end{bmatrix}\\
    			&={\rm det}\begin{bmatrix}
    				Q(G)& 0  &0  &\cdots   & 0 & 0 & 0\\   
    				I_{n} &I_{n}  & 0 &\cdots   &0  & 0 & 0\\
    				0&  I_{n}&  I_{n}&  \cdots & 0 & 0 &0 \\
    				\vdots &\vdots   &\vdots   & \ddots  &\vdots   & \vdots  &\vdots  \\
    				0&  0&  0&  \cdots & I_{n}& 0 &0 \\
    				0& 0 &  0&  \cdots &  I_{n}& I_{n} &0 \\
    				0&  0&  0&  \cdots &  0& I_{n} &I_{n}
    			\end{bmatrix}\\
    		&=\pm 2.
    		\end{aligned}
    	\end{equation*}
    	
    	Now assume that $a_{0}^{t-1} =\pm 2$ $(t\ge 3)$ is true, we will prove that $a_{0}^{t} =\pm 2$ is true.
    	\begin{equation*}
    		\begin{aligned}
    			a_{0}^{(t)}&={\rm det}\begin{bmatrix}
    				Q(G\circ P_{k}^{t-1} )+I_{n}& I_{n}  &0  &\cdots   & 0 & 0 & 0\\   
    				I_{n} &2I_{n}  & I_{n} &\cdots   &0  & 0 & 0\\
    				0&  I_{n}&  2I_{n}&  \cdots & 0 & 0 &0 \\
    				\vdots &\vdots   &\vdots   & \ddots  &\vdots   & \vdots  &\vdots  \\
    				0&  0&  0&  \cdots &  2I_{n}& I_{n} &0 \\
    				0& 0 &  0&  \cdots &  I_{n}& 2I_{n} &I_{n} \\
    				0&  0&  0&  \cdots &  0& I_{n} &I_{n}
    			\end{bmatrix}\\
    			&={\rm det}\begin{bmatrix}
    				Q(G\circ P_{k}^{t-1} )& 0  &0  &\cdots   & 0 & 0 & 0\\   
    				I_{n} &I_{n}  & 0 &\cdots   &0  & 0 & 0\\
    				0&  I_{n}&  I_{n}&  \cdots & 0 & 0 &0 \\
    				\vdots &\vdots   &\vdots   & \ddots  &\vdots   & \vdots  &\vdots  \\
    				0&  0&  0&  \cdots & I_{n}& 0 &0 \\
    				0& 0 &  0&  \cdots &  I_{n}& I_{n} &0 \\
    				0&  0&  0&  \cdots &  0& I_{n} &I_{n}
    			\end{bmatrix}\\
    		&=\pm 2.
    		\end{aligned}
    	\end{equation*}
    	This completes the proof.
    \end{proof}
    
    \paragraph{Theorem 4.6.} Let $G$ be a graph with ${\rm det}W_{Q}(G)=\pm 2^{\frac{3n-2}{2} } $ ($n$ is even) and $a_{0}=\pm2$. Then every graph in the infinite sqeuence $G\circ P_{k}^{t} $ ($k=2,3;~t \ge 1$) are $DGQS$.
    \begin{proof}
    	According to Theorem 2.4 and Lemma 4.5, we can get
    	$${\rm det}W_{Q} (G\circ P_{2}^{t}  )=\pm a_{0}({\rm det}W_{Q} (G\circ P_{2}^{t-1}  )) ^{2} =\cdots =\pm 2^{3n2^{t-1}-1 } ,$$
    	and 
    	$${\rm det}W_{Q} (G\circ P_{3}^{t}  )=\pm (a_{0})^{2}({\rm det}W_{Q} (G\circ P_{2}^{t-1}  )) ^{3} =\cdots =\pm 2^{\frac{3^{t+1}n-2 }{2} } .$$
    	Notice that the number of vertices in $G\circ P_{2}^{t}$ is $2^{t}n$, and the number of vertices in $G\circ P_{3}^{t}$ is  $3^{t}n$. 
    	Thus, we can get that $\frac{{\rm det}W_{Q}(G\circ P_{2}^{t} )}{2^{\frac{3n2^{t}-2 }{2} } }=\pm 1$ and $\frac{{\rm det}W_{Q}(G\circ P_{3}^{t} )}{2^{\frac{3^{t+1}n-2 }{2} } }=\pm 1 $. Therefore, $G\circ P_{2}^{t}$ and $G\circ P_{3}^{t}$ are $DGQS$ from Theorem 4.1.\\
    \end{proof}
    
    Remark that, from Theorem 4.6, it is very easy to construct families of infinite $DGQS$ graphs. Starting from a small $DGQS$ graph $G$ , then we can obtain the following infinite sequence:
    $$G\circ P_{k},~(G\circ P_{k})\circ P_{k} ,\dots ,$$
    where $k=2$ or $k=3$. Clearly, every graph in this sequence is $DGQS$.\\
    
    We conclude this section by posing the following question.
    
     \paragraph{Question 4.7.} Let $G$ be a graph with ${\rm det}W_{Q}(G)=\pm 2^{\frac{3n-2}{2} } $ and $a_{0}=\pm2$. Whether $G\circ P_{k}^{t}$ is $DGQS$ whenever $t\ge 1$ and $k\ge 4$?

    \section{$DGQS$-property of $G\circ P_{2}$ for $G\in \mathcal{F}_{Q,n} $}
    
    In this section, we consider when $G\circ P_{2}$ is $DGQS$ for a graph $G\in \mathcal{F} _{Q,n}$. In the following, we will give some partial answers to this question.
    We first recall some known results in \cite{Qiu2019} that will be used later.
    
    \paragraph{Lemma 5.1} (\cite{Qiu2019}).
    Let $G$ be a graph such that ${\rm det}W_{Q}(G)\ne 0 $. Then there exists a graph $H$ such that $G$ and $H$ are cospectral with respect to the generalized $Q$-spectrum if and only if there exists a rational orthogonal matrix $U$ such that $U^{T}Q(G)U=U(H) $ and $Ue=e$.\\
    
    Define
    $$\Gamma (G)=\left \{ U\in O_{n}(Q)\mid U^{T}Q(G)U=Q(H)~{\rm for~some~graph}~H~{\rm and}~Ue=e  \right \},$$
    where $O_{n}(Q)$ denotes the set of all orthogonal matrices with rational entries.
    
    \paragraph{Lemma 5.2} (\cite{Qiu2019}). Suppose that ${\rm det}W_{Q}(G)\ne0$. Then $G$ is $DGQS$ if and only if $\Gamma (G)$ contains only permutation matrices.
    
    \paragraph{Definition 5.3} (\cite{Qiu2019}). Let $U$ be an orthogonal matrix with rational entries. The level of $U$, expressed as $\ell(U) $ or simply $\ell$, is the smallest positive integer $k$ such that $kU$ is an integral matrix.\\
    
   In order to make the following proof easier, we need to define the modified $Q$-walk matrix $\widetilde{W}_{Q} (G)=\left [ e,~\frac{Qe}{2},~\dots ~,\frac{Q^{n-1}e }{2}   \right ]$ or simply $\widetilde{W}_{Q}$.
    
    \paragraph{Lemma 5.4.}Let the modified $Q$-walk matrices of $G$ and $G\circ P_{2}$ be $\widetilde{W}_{Q}(G)$ and $\widetilde{W}_{Q}(G\circ P_{2})$, respectively. Let $p$ be any prime number. If $\widetilde{W}_{Q}^{T}  (G)\alpha \equiv 0~({\rm mod}~p)$, then
    $$\widetilde{W}_{Q}^{T} (G\circ P_{2} )  \begin{bmatrix}
    	\alpha \\
    	0
    \end{bmatrix}\equiv\widetilde{W}_{Q}^{T} (G\circ P_{2} )  \begin{bmatrix}
    	0\\
    	\alpha 
    \end{bmatrix}\equiv 0~({\rm mod}~p).$$
    \begin{proof}
    	We can first see that, for any integer $k\ge 0$,
    	$$\frac{1}{2} Q^{k} (G\circ P_{2} )=\begin{bmatrix}
    		Q(G)+I_{n} & I_{n}\\
    		I_{n}&I_{n}
    	\end{bmatrix}^{k} =\begin{bmatrix}
    		f(Q)&g(Q) \\
    		g(Q)&h(Q)
    	\end{bmatrix} ,$$
    	where polynomials $f(x)$, $g(x)$ and $h(x)$ have degrees of $k$, $k-1$ and $k-2$, respectively.
    	Due to $\widetilde{W}_{Q}^{T}  (G)\alpha \equiv 0~{(\rm mod~p)}$, then we can get 
    	\begin{equation}\label{9}
    		e^{T} \alpha \equiv 0,~e^{T}\frac{Q^{k} }{2}\alpha \equiv 0~({\rm mod}~p).
    	\end{equation}
    	It follows that 
    	\begin{equation*}
    		\begin{aligned}
    			\left [ e^{T},e^{T}   \right ] \frac{1}{2} Q^{k} (G\circ P_{2} )\begin{bmatrix}
    				\alpha \\
    				0
    			\end{bmatrix}&=\left [ e^{T},e^{T}   \right ]\frac{1}{2}\begin{bmatrix}
    				f(Q)&g(Q) \\
    				g(Q)&h(Q)
    			\end{bmatrix}\begin{bmatrix}
    				\alpha \\
    				0
    			\end{bmatrix}\\
    			&=e^{T} \frac{1}{2} f(Q)\alpha +e^{T} \frac{1}{2} g(Q)\alpha\\
    			&\equiv 0~({\rm mod}~p).
    		\end{aligned}
    	\end{equation*}
    	In the same way, we have $\left [ e^{T},e^{T}   \right ] \frac{1}{2} Q^{k} (G\circ P_{2} )\begin{bmatrix}
    		0 \\
    		\alpha
    	\end{bmatrix}\equiv 0~({\rm mod}~p).$ This completes the proof.
    \end{proof}
    
    \paragraph{Theorem 5.5.} Given a graph $G\in \mathcal{F} _{Q,n}$ ($n$ is even) and $a_{0}=\pm 2$. Let $p$ be any odd prime factor of ${\rm det}\widetilde{W} _{Q} (G)$. Suppose that both of the following conditions are true:
    \begin{enumerate}[(i)]
    	\item $\widetilde{W}^{T} _{Q}(G)\alpha \equiv 0({\rm mod}~p)$ and $\alpha ^{T}\alpha\not \equiv  0~({\rm mod}~p) $.
    	\item  $p\nmid {\rm Res}(P_{Q}(x),\frac{x-1}{2} )$, where $P_{Q}(x)$ is characteristic polynomial of $Q(G)$.
    \end{enumerate}
    Then the rooted product graph $G\circ P_{2}$ is $DGQS$.
    \begin{proof}
         Remark that the proof of this theorem is similar to Theorem 3.9 in \cite{Mao2022}. It follows from Theorem 2.4 that ${\rm det}W_{Q}(G\circ P_{2} )=\pm 2^{3n-2}p_{1}^{2}p_{2}^{2}\cdots p_{s}^{2}$. By the definition of the modified $Q$-walk matrix, we have ${\rm det}\widetilde{W}_{Q}(G\circ P_{2} )=\pm 2^{n}p_{1}^{2}p_{2}^{2}\cdots p_{s}^{2}$.
         
         First, the case of $p=2$ can easily be ruled out from Theorem 4.1. So we just need to rule out the odd prime factor $p_{1},p_{2},\dots ,p_{s}$. Let $p>2$ be any prime factor of ${\rm det}W_{Q}(G\circ P_{2} )$. Suppose that $\widehat{U}$ is
         the rational orthogonal matrix with level $l$ such that $\widehat{U} ^{T} Q(G\circ P_{2} )\widehat{U}=A(\Gamma )$, where $\Gamma$ is a graph being generalized $Q$-cospectral with $G\circ P_{2}$. We will prove $p\nmid \ell$. By contradiction, assume that $p\mid \ell$. Let $\nu $ be any column of $\ell\widehat{U} $. Then we obtain
         \begin{equation*}
         	\widetilde{W} _{Q}(G\circ P_{2} )\nu \equiv 0,~\nu ^{T}\nu =0~({\rm mod}~p).
         \end{equation*}
         Observe that ${\rm rank}_{p}\widetilde{W}_{Q}(G)=n-1$. It follows that the nullspace of $\widetilde{W}_{Q}^{T}(G)$ is $1$-dimensional.
         Let $0\ne \alpha \in {\rm Null}(\widetilde{W}_{Q}^{T}(G))$. Then from Lemma 5.3, we know that $\nu$ can be represented as a linear combination of $\begin{bmatrix}
         	\alpha \\
         	0 
         \end{bmatrix}$ and $\begin{bmatrix}
         0\\
         \alpha 
     \end{bmatrix} $. There are some $k$ and $l$ such that $\nu =\begin{bmatrix}
     k\alpha \\
     l\alpha 
     \end{bmatrix} $, where $k,l\in \mathbb{F}_{p}$. Then we have 
     \begin{equation*}
     	\ell\widehat{U}=\begin{pmatrix}
     		k_{1}\alpha  &k_{2}\alpha    &\cdots   &k_{2n}\alpha   \\
     		l_{1} &l_{2}   &\cdots   &l_{2n}\alpha  
     	\end{pmatrix} ~({\rm mod}~p).
     \end{equation*}
      Observe that $\widehat{U} $ is an orthogonal matrix. Thus, $(k_{i}^{2}+l_{i}^{2}   )\alpha ^{T}\alpha \equiv 0~({\rm mod}~p) $ and $(k_{i}k_{j} +l_{i}l_{j}   )\alpha ^{T}\alpha \equiv 0~({\rm mod}~p) $. Since $\alpha ^{T}\alpha\not \equiv  0~({\rm mod}~p) $, then $k_{i}^{2}+l_{i}^{2}\equiv 0~({\rm mod}~p) $. Obviously, there exists some $k_{i},~l_{i}$ with $k_{i}\not \equiv  0,~l_{i}\not \equiv0~({\rm mod}~p)$. If $p\equiv 3~({\rm mod}~4)$, then $(k_{i}l_{i}^{-1}  )^{2} +1\equiv 0~({\rm mod}~p)$, this is contradictory to the fact that $k_{i}l_{i}^{-1}$ is a rational number. Now we may suppose that $p\equiv 1~({\rm mod}~4)$, then the congruence equation $x^{2} +1\equiv 0~({\rm mod}~p)$. Hence, ${\rm rank }_{p}(\ell \widehat{U}) =1$. Suppose the number of non-zero elements in $\alpha$ is $r$. Then we know that $\ell \widehat{U}$ can be written as
      \begin{equation*}
      	\ell \widehat{U}\equiv \begin{pmatrix}
      		k\alpha & m_{1}k\alpha   &\cdots   &0  & \cdots  & 0\\
      		l\alpha & m_{1}l\alpha   &\cdots   & 0 & \cdots  &0
      	\end{pmatrix}~({\rm mod}~p).
      \end{equation*}
      It follows from $\widehat{U} ^{T} Q(G\circ P_{2} )\widehat{U}=A(\Gamma )$ that
      \begin{equation*}
          \begin{pmatrix}
          	Q(G)+I& I\\I
          	&I
          \end{pmatrix}\equiv b_{0}\begin{pmatrix}
          	k\alpha \\
          	l\alpha 
          \end{pmatrix} +b_{1} \begin{pmatrix}
          	m_{1} k\alpha \\
          	m_{1} l\alpha 
          \end{pmatrix} +\cdots +b_{2r-1} \begin{pmatrix}
          	m_{2r-1} k\alpha \\
          	m_{2r-1}l\alpha 
          \end{pmatrix} .
      \end{equation*}  
         Thus,
         \begin{equation*}
         	kQ(G)\alpha +l\alpha =\lambda k\alpha ,~k\alpha +l\alpha =\lambda l\alpha ,
         \end{equation*}
     where $\lambda :=b_{0}+m_{1} b_{1}+\cdots +m_{2r-1}b_{2r-1}.$ The first equation implies that $Q(G)\alpha +k^{-1}l\alpha \equiv \lambda \alpha \equiv l^{-1}k\alpha +\alpha~({\rm mod}~p) $. The second equation implies that $k+l\equiv \lambda l~({\rm mod}~p)$.  That is, $Q(G)\alpha=(l^{-1}k-k^{-1}l+1  )\alpha  ~({\rm mod}~p)$. Now let $\lambda =kl^{-1} $. Then we obtain 
     \begin{equation}
     	\lambda ^{2} +1\equiv 0~({\rm mod}~p)~{\rm and}~Q(G)\alpha =(2\lambda +1) \alpha .
     \end{equation}
     It is clear that the Equation $(11)$ is equivalent to the following condition: $x^{2}+1$ and $P_{Q}(\frac{x-1}{2}) $ has a common root over $\mathbb{F}_{p}$, i,e. $p\mid {\rm Res}(P_{Q}(x), \frac{x-1}{2})$. This contradicts the condition (ii). Therefore, $p\nmid \ell $. Then $\ell=1$ and $\widehat{U} $ is a permutation matrix.
     This completes the proof.
    \end{proof}
    
    \paragraph{Corollary 5.6.}Given a graph $G\in \mathcal{F} _{Q,n}$ ($n$ is even) and ${\rm det}\widetilde{W} _{Q}(G)=\pm 2p_{1}p_{2} \cdots p_{t}$ ($n$ is even and $p_{1},p_{2},\dots,p_{t}$ are distinct odd primes) and $a_{0}=\pm 2$. If $p_{i}  \equiv 3~({\rm mod} ~4) $ for $i=1,2,\dots,t$, then the graph $G\circ P_{2}$ is $DGQS$.\\

    \noindent
    \textbf{Declaration of competing interest} The authors report no potential conflicts of interest.\\
    \\
    \textbf{Acknowledgements} This work was in part supported by the National Natural Science Foundation of China (Nos. 11801521, 12071048).

\end{document}